\documentclass[11pt,leqno]{article}
\usepackage{amsmath}
\usepackage{amsthm}

\numberwithin{equation}{section}
\let\sect=\section

\newtheorem{theorem}{Theorem}[section]
\newtheorem{corollary}[theorem]{Corollary}
\newtheorem{lemma}[theorem]{Lemma}
\newtheorem{proposition}[theorem]{Proposition}
\newtheorem{claim}[theorem]{Claim}
\newtheorem{example}[theorem]{\sl Example}

\theoremstyle{definition}
\newtheorem{remark}[theorem]{Remark}

\newcommand{\lf}{\left\lfloor}

\newcommand{\rf}{\right\rfloor}

\newcommand{\EE}{{\bf  E}}

\newcommand{\RR}{{\bf  R}}

\newcommand{\ZZ}{{\bf  Z}}

\newcommand{\UU}{{\bf  U}}
\newcommand{\XX}{{\bf  X}}
\newcommand{\YY}{{\bf  Y}}

\newcommand{\uu}{{\bf  u}}
\newcommand{\xx}{{\bf  x}}
\newcommand{\yy}{{\bf  y}}

\newcommand{\ww}{{\bf  w}}

\newcommand{\Var}{{\bf Var}}

\newcommand{\id}{{\rm id}}
\newcommand{\rev}{{\rm rev}}

\newcommand{\Kt}{{\tilde{K}}}

\newcommand{\SSS}{{\cal S}} 

\newcommand{\Lc}{{\cal L}}

\newcommand{\pc}{{p_{\mbox{\rm \scriptsize CFTP}}}}
\newcommand{\pf}{{p_{\mbox{\rm \scriptsize FMMR}}}}
\newcommand{\Tc}{{T_{\mbox{\rm \scriptsize CFTP}}}}
\newcommand{\Tf}{{T_{\mbox{\rm \scriptsize FMMR}}}}
\newcommand{\Zc}{{Z_{\mbox{\rm \scriptsize CFTP}}}}
\newcommand{\Zf}{{Z_{\mbox{\rm \scriptsize FMMR}}}}

\newcommand{\Uc}{{\cal U}}
\newcommand{\Xc}{{\cal X}}

\newcommand{\zh}{\hat{0}}
\newcommand{\oh}{\hat{1}}

\newcommand{\begp}{\begin{proposition}}
\newcommand{\enp}{\end{proposition}}
\newcommand{\begt}{\begin{theorem}}
\newcommand{\ent}{\end{theorem}}
\newcommand{\begl}{\begin{lemma}}
\newcommand{\enl}{\end{lemma}}
\newcommand{\begc}{\begin{corollary}}
\newcommand{\enc}{\end{corollary}}
\newcommand{\begcl}{\begin{claim}}
\newcommand{\encl}{\end{claim}}
\newcommand{\begr}{\begin{remark}}
\newcommand{\enr}{\end{remark}}
\newcommand{\begal}{\begin{algorithm}}
\newcommand{\enal}{\end{algorithm}}
\newcommand{\begd}{\begin{definition}}
\newcommand{\enf}{\end{definition}}
\newcommand{\begx}{\begin{example}}
\newcommand{\enx}{\end{example}}
\newcommand{\bega}{\begin{array}}
\newcommand{\ena}{\end{array}}

\def\rompar(#1){\textup(#1\textup)}    

\newcommand{\refS}[1]{Section~\ref{#1}}
\newcommand{\refT}[1]{Theorem~\ref{#1}}

\newcommand\ie{i.e.\spacefactor=1000}

\begin{document}

\setcounter{page}{0}
\thispagestyle{empty}

\begin{center}
{\Large \bf Speeding up the FMMR perfect sampling algorithm: \\ A case
study revisited\\ }
\vspace{.15in}
\normalsize

\vspace{4ex}
{\sc Robert P.\ Dobrow} \\
\vspace{.1in}
Mathematics and Computer Science Department\\
\vspace{.1in}
Carleton College \\
\vspace{.1in}
{\tt rdobrow@carleton.edu}
and \\
{\tt http://www.mathcs.carleton.edu/faculty/bdobrow/} \\
\vspace{.2in}
{\sc and} \\
\vspace{.2in}

{\sc James Allen Fill\footnotemark} \\
\vspace{.1in}
Department of Mathematical Sciences \\
\vspace{.1in}
The Johns Hopkins University \\
\vspace{.1in}
{\tt jimfill@jhu.edu} and {\tt http://www.mts.jhu.edu/{\~{}}fill/} \\
\end{center}
\vspace{4ex}

\begin{center}
{\sl ABSTRACT} \\
\end{center}
\vspace{.05in}

\begin{small}
In a previous paper by the second author,
two Markov chain Monte Carlo perfect sampling algorithms---one called
coupling from the
past (CFTP) and the other (FMMR) based on rejection sampling---are compared
using as a
case study the move-to-front (MTF) self-organizing list chain.  Here
we
revisit that case
study and, in particular, exploit the dependence of FMMR on the
user-chosen
initial state.
We give a stochastic monotonicity result for the running time of FMMR
applied to MTF and
thus identify the initial state that gives the stochastically smallest
running time; by
contrast, the initial state used in the previous study gives the
stochastically {\em
largest\/} running time.  By changing from worst choice to best choice of
initial state we
achieve remarkable speedup of FMMR for MTF; for example, we
reduce the running
time (as measured in Markov chain steps)
from exponential in the length~$n$ of the list nearly down to~$n$
when the
items in the
list are requested according to a geometric distribution.  For this same
example, the
running time for CFTP grows exponentially in~$n$.
\medskip
\par\noindent
{\em AMS\/}
2000 {\em subject classifications.\/}  Primary 60J10, 68U20; secondary 60G40,
68P05, 68P10, 65C05, 65C10, 65C40.
\medskip
\par\noindent
\emph{Date.\/} May~14, 2002.
\end{small}

\footnotetext[1]{Research for this author supported by NSF grants
DMS--9803780 and DMS--0104167,
and for both authors by The Johns Hopkins University's Acheson J.~Duncan Fund for the
Advancement of Research in Statistics.}

\newpage
\addtolength{\topmargin}{+0.5in}

\sect{Introduction and summary}
\label{S:intro}

Perfect sampling has had a substantial impact on
the world
of Markov Chain Monte Carlo (MCMC).  In MCMC, one is interested in
obtaining a sample from a
distribution~$\pi$ from which it is computationally difficult (or even
infeasible) to
simulate directly.  One constructs a Markov chain whose stationary
distribution is~$\pi$ and
after running the chain ``a long time'' takes an outcome from the chain as
an (approximate)
observation from~$\pi$.  Propp and Wilson~\cite{propp-wilson} (see
also \cite{PWuser}~\cite{PWhowto}~\cite{wilson})
and Fill~\cite{fill-ps1} have devised algorithms to use Markov
chain transitions to produce observations {\em exactly\/} from~$\pi$,
without {\em a
priori\/} estimates of the mixing time of the chain; the applicability of
the latter
algorithm has recently been extended by Fill, Machida, Murdoch, and
Rosenthal~\cite{FMMR},
and so we will use the terminology ``FMMR algorithm.''  Although the
two algorithms are
based on different ideas---Propp and Wilson use coupling from the past
(CFTP) while FMMR is
based on rejection sampling---there is a simple connection between the two,
discovered
in~\cite{FMMR} and reviewed below.  For further general discussion of
perfect sampling using
Markov chains, consult the annotated bibliography maintained on the Web by
Wilson~\cite{wilsonbib}.

Much of the discussion comparing the two algorithms has focused on the issue of
``interruptibility.''  FMMR has the feature that the output and the
running time---when
measured in number of Markov chain steps---are independent random
variables.  Thus, for
instance, an impatient user who interrupts a run of the algorithm after any
fixed number
of steps and restarts the procedure does not introduce any bias into the
output.  This is
not so for CFTP.  On the other hand, for many practical applications CFTP
is considerably
easier to implement, since (see~\cite{FMMR}) FMMR requires the user to be
able (i)~to
generate a trajectory from the time-reversal of the basic chain, and
(ii)~to build couplings
``{\em ex post facto\/},'' i.e.,\ to perform certain imputation steps; CFTP
requires neither
ability.

\begin{remark}
There {\em is\/} a need for time-reversal generation (of an auxiliary
chain) and for {\em
ex post facto\/} coupling in an extension of CFTP known as coupling into
and from the
past, introduced (under a different name) by Kendall~\cite{kendall}.   (See
also Section~1.9.3 in~\cite{wilson}.)
\end{remark}

In this paper we focus on the running time of the two algorithms (but  the
non-interruptibility of CFTP will turn out to play a key role).  In
previous case-study
comparisons \cite{fill-ps1}~\cite{fill-mtf-ps}, the running times (and
memory requirements)
have been found to be not hugely different, but CFTP has had the edge.  In
this paper, by
revisiting the case study of~\cite{fill-mtf-ps}, we show that, at least in
some cases, FMMR
can be made to have much smaller running time than CFTP.

The general observation that we exploit---one very closely related to
Remark~6.9(c) and
Section~8.2 of~\cite{FMMR}---is the following.  Given a target
distribution~$\pi$, let~$\pc$
denote the probability that CFTP terminates successfully (coalesces)
over a
fixed time
window (and outputs a sample from~$\pi$).  Similarly, let~$\pc(z)$ denote
the conditional
probability of coalescence over the time window, given that the state
(call it~$\Zc$)
ultimately output by CFTP (after extending the time window into the
indefinite past) is~$z$.
Let~$\pf(z)$ denote the conditional probability that FMMR terminates
successfully over the
same time window, given that the initial state (call it~$\Zf$) is~$z$.
Then, as we show in
Theorem~\ref{runtime}, $\pc(z) \equiv \pf(z)$.
That is (now letting the time window vary),
if~$\Tc$ and~$\Tf$ denote the respective running times of CFTP and FMMR,
then conditional
running time distributions agree:
$$
\Lc(\Tc\,|\,\Zc = z) \equiv \Lc(\Tf\,|\,\Zf = z),
$$
where $\Lc(X)$ denotes the distribution (law) of the random variable
$X$.
As a consequence, $\pc = \EE_{\pi}\,[\pf(\Zf)]$; that is,\ $\Lc(\Tc)$ is the
$\pi$-mixture of
the distributions $\Lc(\Tf\,|\,\Zf = z)$.

The important point here is that, except in the rare instance that CFTP is
interruptible for
the chain of interest (i.e.,\ that~$\Tc$ and~$\Zc$ are independent),
for at least one time
window there must exist at least one initial state~$z$ for which $\pf(z) >
\pc$. 

The move-to-front (MTF) process is a nonreversible Markov chain on the permutation
group~${\SSS}_n$.
The two algorithms have been compared for MTF in a previous
paper~\cite{fill-mtf-ps}.
In that paper, the initial state for FMMR was taken to be the identity
permutation, and it
was then found, roughly speaking (see Table~1 and Section~5 therein),
that~$\Tc$ and~$\Tf$
are of the same size.  In this paper, we will revisit that case study  and
establish a stochastic monotonicity result for $\Lc(\Tc\,|\,\Zc = z)$ in~$z$.
It turns out, in
particular, that the identity permutation is the {\em worst\/} choice  of
initial state!
When we choose instead the reversal permutation, which is the best choice,
we obtain a
(sometimes huge) speedup for FMMR.  (See Table~1, which will be explained
more fully in
Section~\ref{speedup}.
Notice that for geometric weights, the change in starting state
reduces~$\Tf$ from exponential in~$n$ to about~$n$.)
The gains obtained by using the optimal~$z$ are sufficiently dramatic
that, when measured
in Markov chain steps, the resulting {\em worst-case\/} running time
for FMMR (worst over
choice of request weights) equals the {\em best-case\/} running time
for CFTP:\
see Remark~\ref{extremes}(b).

We temper our enthusiasm, however, by recognizing
that our MTF example is somewhat artificial on two counts.  Firstly, as
discussed in the introduction to~\cite{fill-mtf-ps}, there are algorithms for
sampling from the MTF stationary distribution which are both more
elementary (in particular,
not involving Markov chains) and more efficient than either CFTP or FMMR.
So we do not
recommend applying either CFTP or FMMR to MTF in practice. Our goal here
is to illustrate how
judicious choice of starting state for FMMR can greatly improve its performance.

Secondly, MTF has the (evidently rare) property that one can obtain an
exact analysis of the running time distribution for FMMR for {\em every\/}
choice of initial state~$z$.  We do not yet know whether our speedup
ideas help in any real applications.
We hope, however, that the ideas in this paper will stimulate further  research
on FMMR by pointing to the possibility of speedup of the algorithm.

We briefly review the two perfect sampling algorithms and their general
connection in
Section~\ref{psreview}.
The move-to-front rule is reviewed in Section~\ref{mtf}.  Our new results
are given in
Section~\ref{speedup}.  A somewhat different approach to speeding up FMMR is given
in Section~\ref{last}.

\sect{Perfect sampling}\label{psreview}
We briefly review the CFTP and FMMR algorithms (omitting a few of the  finer
measure-theoretic details, which are irrelevant anyway for finite-state chains).
We assume
that our Markov chain
$\XX$ can be written in the {\em stochastic recursive sequence\/} form
\begin{equation}
\label{recursive}
\XX_s = \phi(\XX_{s-1}, \UU_s),
\end{equation}
where~$\phi$ is called the {\em transition rule\/} and $(\UU_s)$ is an i.i.d.\ sequence.
We further assume that our Markov chain has 
finite state space~$\Xc$ and is irreducible and aperiodic with (unique)
stationary distribution $\pi$.

\subsection{CFTP}

For a fixed positive integer~$t$, and a Markov chain with~$n$ states,
start~$n$ copies of the chain at time $- t$ from each of the~$n$ states,
coupling the transitions by means of the transition rule~$\phi$,
and running the chains until time~$0$.  If all copies of the chain agree
at time~$0$, we say that the trajectories have {\em coalesced\/} and
return the common value, say~$\ZZ$.  If the chains do not agree,
then increment~$t$ and
restart the procedure, using for
common values of~$s$ the same values
of~$\UU_s$ used in the previous step; again,
check for coalescence.
The {\em running time\/} of
the algorithm we
define to be the smallest integer~$t$ for which coalescence occurs.
If we assume the algorithm terminates with probability~$1$,
then $\ZZ \sim \pi$ exactly.

There is a rich source of papers, primers, and applications of CFTP.
The best initial reference is the ``Perfectly random sampling
with Markov chains'' Web site
maintained by David Wilson at
{\tt http://www.dbwilson.com/exact/}.

\subsection{FMMR}\label{subfmmr}
Given a Markov chain with transition matrix $K$, recall that the
{\em time-reversal\/} chain has transition matrix
$\tilde{K}$ which
satisfies
$$ \pi(x) K(x,y) = \pi(y) \tilde{K}(y,x) \mbox{ for all } x, y .$$

The FMMR algorithm has two stages: First, choose an initial state $\XX_0$.
Run the time-reversed chain $\tilde{K}$, obtaining $\XX_0, \XX_{-1},
\ldots$ in succession. Then (conditionally given the $\XX$-values)
generate $\UU_0, \UU_{-1}, \ldots$ independently, with $\UU_s$ chosen
from its conditional distribution given (\ref{recursive}) for
$s=0, -1, \ldots$.  (One says that the values $\UU_s$ are \emph{imputed}.)\ \ For $t=0, 1,
\ldots$, and for each state $x$ in the state space, set $\YY_{-t}^{(-t)}(x) := x$ and,
inductively,
$$ \YY_s^{(-t)}(x) := \phi ( \YY_{s-1}^{(-t)}(x), \UU_s), ~~~ -t +1
\leq s \leq 0.$$
We will sometimes refer to the realization of the chain~$\XX$ as the
\emph{backward} trajectory, and to the  realizations of the chains $\YY(x)$ as the 
{\em forward} trajectories.
The {\em running time\/} of the algorithm we define to be the smallest $t^*$ such
that $\YY_0^{(-t^*)}(x)$ agree for every~$x$ in the state space (and
hence all equal $\XX_0$).  In this case the algorithm reports $\XX_{-t^*}$
as an observation from $\pi$.

\begin{remark}
The algorithms are presented here in their most general,
``vanilla'' versions.  A large amount of research has gone into improving 
both algorithms and tailoring them for specific applications.  In
particular, to improve performance a ``doubling trick'' is suggested for
both algorithms whereby instead of incrementing $t$ by one at each step,
$t$ is successively doubled.
Since this affects the number of Markov chain steps taken only by
constant factors, we shall for our theoretical analysis stick to the
``vanilla'' versions.
\end{remark}

\begin{remark}
For most chains of interest, the state space is very large and the
implementations presented here (running copies of the chain from every
state in the state space) are not feasible.  However, for a large class
of cases where a form of monotonicity holds, the algorithms become
practical.

Given a Markov chain with transition matrix $K$, we say that we are
in the
{\em (realizably) monotone case\/} if the following conditions hold.
The state space is
a partially
ordered set
$(\Xc, \leq)$.  There exist (necessarily unique) minimum and maximum elements
in the state
space, denoted $\zh$ and $\oh$, respectively. There exists a {\em monotone
transition rule\/}
$\phi$ for the chain.  Such a rule is a function
$\phi: \Xc \times \Uc
\to \Xc$
together with a random variable $\UU$ taking values in a probability
space
$\Uc$ such that
(i) $\phi(\xx, \uu) \leq \phi(\yy,\uu)$ for all $\uu \in \Uc$ whenever
$\xx \leq \yy$; and (ii)~for
each $\xx \in \Xc$, $P(\phi(\xx,\UU ) = \yy) = K(\xx,\yy)$ for all $\yy \in
\Xc$.

When in the monotone case, for CFTP
one only needs to follow two trajectories
of the chain, one started at time $-t$ from~$\zh$ and the other from~$\oh$,
since all other trajectories are sandwiched between these.
Likewise, in
the second phase of FMMR, one only needs to run the $\YY$-chain
from states $\zh$ and $\oh$.
\end{remark}

Although the two algorithms are based on different conceptual underpinnings,
our first theorem  highlights an important connection between them.
Roughly, the distribution of the running time for CFTP is equal
to the
stationary mixture,
over initial states, of the distributions of the running time for FMMR.
This is given as Remark~6.9(c) in~\cite{FMMR}, but we wish to emphasize
its importance and so recast it as a
theorem. We recall our notation from \refS{S:intro}.
For a fixed time window, let $\pc(z)$ denote the probability that CFTP
coalesces given that the
state (call it $\Zc$) ultimately output by CFTP is~$z$,
and let $\pc$ denote the corresponding {\em unconditional\/} probability.
Let $\pf(z)$ denote the conditional
probability that FMMR coalesces given that the initial state (call it  $\Zf$) is $z$.
Let $\Tc$ and $\Tf$ denote the respective running times of CFTP and FMMR
(now letting the time window vary).

\begin{theorem}\label{runtime}
We have 
\begin{equation}\label{run1}
\pc(z) = \pf(z) \mbox{\rm \ for $\pi$-almost every~$z$};
\end{equation}
\begin{equation}\label{run3}
\pc = {\bf E}_{\pi}[\pf(\Zf)];
\end{equation}
\begin{equation}\label{run2}
\Lc(\Tc | \Zc = z) = \Lc(\Tf | \Zf = z) \mbox{\rm \ for $\pi$-almost every~$z$};
\end{equation}
and
\begin{equation}\label{mixture}
\mbox{\rm $\Lc(\Tc)$ is the $\pi$-mixture (over~$z$) of $\Lc(\Tf\,|\,\Zf = z)$.}
\end{equation}
\end{theorem}

The result holds in the most general setting, not restricted either to
finite-state chains or to
monotone transition rules.  It is a consequence of the discussion in
Sections~6.2 and~8.2 of
\cite{FMMR}.  For the reader's convenience we give here a simple
proof for the discrete case.

\begin{proof}
Each
iteration of FMMR is an implementation of rejection sampling  (see,\ 
e.g.,\ Devroye~\cite{dev} for background).  The goal is to use an observation from
$\Kt^t(z,\cdot)$ to
simulate one from $\pi$.  One obtains an upper bound~$c$ on
$\max_x \pi(x) / \Kt^t(z, x)$,
generates~$x$ with probability $\Kt^t(z,x)$, and accepts $x$ as an
observation from $\pi$ with
probability $c^{-1} \pi(x)/\Kt^t(z,x)$.  The unconditional probability
of acceptance is then
$1/c$.  Observe that, for every~$x$,
$$  \frac{\pi(x)}{\Kt^t(z,x)} = \frac{\pi(z)}{K^t(x,z)}
  \leq \frac{\pi(z)}{P(\mbox{coalescence to~$z$})} := c,$$
where ``coalescence to~$z$'' refers, of course, to coalescence over
the given time window of length~$t$.
Thus for the desired conditional acceptance probability given~$x$ we
can use
$$
\frac{P(\mbox{coalescence to~$z$})}{K^t(x, z)}
  = P(\mbox{coalescence to~$z$}\,|\,\mbox{trajectory from~$x$ ends at~$z$}),
$$
and the FMMR algorithm is designed precisely to implement this.
Thus $\pf(z) = 1 / c$ and hence
\begin{equation}\label{run4}
\pf(z) = \frac{P(\mbox{coalescence to~$z$})}{\pi(z)} = \pc(z).
\end{equation}
Thus, (\ref{run1}) is immediate.
Taking expectations with respect to~$\pi$ gives~(\ref{run3}).
And~(\ref{run2})
[from which~\eqref{mixture} is immediate]
follows from (\ref{run1}) since, for a fixed time window of length~$t$,
$\pf(z)$ [respectively, $\pc(z)$] is the value at~$t$ of the
conditional distribution function
of~$\Tf$ given that $\Zf = z$ (respectively, of~$\Tc$ given that $\Zc = z$).
\end{proof}

\begin{corollary}
If~$\Tc$ and~$\Zc$ are not independent random variables,
then there exist  at least one time window and at least one initial
state~$z$ for which  $\pf(z) > \pc$. 
\end{corollary}

\begin{proof}
This is immediate from~(\ref{run3}).
\end{proof}

The following simple examples are  artificial, but they give a first
demonstration that
judicious choice of starting state can lead to dramatic speedup.  First,
consider a three-state Markov chain with states labeled~$0$, $1$, and~$2$.  Let
$$ K = \left[ \begin{array}{ccc}\epsilon & (1-\epsilon)/2 & (1-\epsilon)/2 \\
\epsilon & 1-\epsilon & 0 \\ \epsilon & 0 & 1- \epsilon \end{array} \right] ,$$
where $\epsilon > 0$ is small.  One checks that this chain is reversible
(but not monotone).
Let $\UU = 0, 1, 2$ with respective probabilities $\epsilon, (1 - \epsilon) / 2,
(1 - \epsilon) / 2$ and use the natural transition rule
$$
\phi(x, 0) = 0\mbox{\ for all~$x$},\ \ \phi(0, 1) = \phi(1, 1)
  = \phi(1, 2) = 1,\ \ \phi(0, 2) = \phi(2, 1) = \phi(2, 2) = 2.
$$
Coalescence occurs over a given time window of length~$t$ if and only
if the value of some~$\UU_s$ in
that window is~$0$; thus $\pc = 1 - (1 - \epsilon)^t$, which
requires~$t$ of order~$1 / \epsilon$ to
become nonnegligible.  On the other hand, 
if FMMR is started in state~$0$, then with
high probability
($\ = 1 - \epsilon$) we'll see (going backward in time) one of the
transitions $1 \leftarrow 0$ or $2 \leftarrow 0$.  If we do, then 
(whichever we see) in the forward phase we impute $\UU = 0$
and hence get coalescence (to state~$0$) 
in one step.

For our second example, consider a Gibbs sampler on an
attractive spin system with~$n$ sites
arranged in a row and left-to-right site-update sweeps.
(Consult, e.g.,\ \cite{fill-ps1}
or~\cite{liggett} for background on attractive
spin systems.)  This gives a monotone, nonreversible chain where~$\zh$
is the state consisting of all $-$'s and $\oh$ is the state of all $+$'s.
Suppose that the Gibbs distribution is such
that there is (i)~a strong external field for spin~$+$
at sites~$1$ through $n - 1$,
(ii)~a much stronger effect of attractiveness throughout the system, and
(iii)~a very much stronger yet
external field for spin~$+$ at site~$n$ (the rightmost
site).  First consider CFTP.  The state~$\oh$ is a state of very high probability
and so the chain won't
budge out of that state for a long time.  On the other hand, from~$\zh$,
in one sweep (a full left-to-right update), we obtain
[with high probability, because of~(ii) and~(iii)] $- - \cdots - - +$.
In the next sweep we obtain [with high probability, because of~(ii) and~(i)]
$- - \cdots - + +$. Continuing, in about $n$ sweeps we obtain $+ + \cdots + + +$;
that is, with high probability we obtain coalescence in $n$ sweeps.
On the other hand, consider FMMR started in~$\zh$.  Here, the reversed
chain is Gibbs sampling with right-to-left updates.  The reversed chain,
started in $\zh$, [with high probability, because of~(iii) and then~(ii)]
flips each site from~$-$ to~$+$ as it moves from right to left.
Hence, we obtain $+ + \cdots + + +$; that is, with high
probability there is coalescence {\em in one sweep\/}.

(Of course, were we to use right-to-left sweeps or reversible sweeps as the sampler,
the relative disadvantage of CFTP would disappear.) 

\begin{remark}\label{nogeneral}
In general, we know of no simpler expression for~$\pf(z)$  than the
ratio in~(\ref{run4}).
In the monotone case, however, when $z = \zh$ or $z = \oh$ we obtain
significant simplification. Indeed, then
$$ \pf(\zh) = \frac{K^t(\oh, \zh)}{\pi(\zh)} = \min_z \frac{K^t(z, \zh)}{\pi(\zh)}
~~~ \mbox{ and } ~~~
\pf(\oh) = \frac{K^t(\zh,\oh)}{\pi(\oh)} = \min_z \frac{K^t(z,
  \oh)}{\pi(\oh)}.$$
Recall that for a Markov chain with transition matrix~$K$ and stationary
distribution~$\pi$,
the {\em separation\/} at time~$t$ given that the chain starts in state~$x$ is
$$ \mbox{sep}_x(t) := 1 - \min_z \frac{K^t(x,z)}{\pi(z)}.$$
Thus, $\pf(z) = 1 - \widetilde{\mbox{sep}}_z (t)$ for $z = \zh, \oh,$
where $\widetilde{\mbox{sep}}$ refers to separation for the transition matrix~$\Kt$.
See (e.g.)~\cite{ald-fill} for more on separation.
\end{remark}

\sect{Move-to-front}\label{mtf}
Let
$(w_1, \ldots, w_n)$ be a probability mass function on $\{1, \ldots, n\}$
with $w_i > 0$ for each~$i$.
In this study we are concerned with generating an observation from the distribution
\begin{equation}
\label{target}
  \pi(z) := \prod_{r = 1}^n \frac{w_{z_r}}{\sum_{j = r}^n w_{z_j}},
                   \ \ z \in {\SSS}_n,
\end{equation}
where~${\SSS}_n$ is the group of permutations of $\{1, \ldots, n \}$.
Consider sampling without replacement from a population of $n$ items,
where item $i$ has
probability $w_i$ of being chosen, $1 \leq i \leq n$. Then the
probability of drawing the~$n$ items in the order~$z$ is given
by~(\ref{target}).

This distribution arises as the limiting distribution of
the much-studied
move-to-front (MTF) process.
The MTF heuristic is used to ``self-organize'' a linear list of data
records in a computer
file.  Let $\{1, \ldots, n\}$ be a set of
records (or rather the set of {\em keys\/}, or identifying labels for
the records), where
record~$i$ has probability $w_i$ of being requested.  At discrete
units of time, and independent of
past requests, item $i$ is requested (with probability $w_i$) and
brought to the front of the list,
leaving the relative
order of the other records unchanged.
The successive orders of the list of records forms an ergodic Markov
chain on the permutation
group ${\SSS}_n$ with stationary distribution~$\pi$.

We will assume that the records have
been labeled so
that $w_1 \geq \cdots \geq w_n > 0$ and refer to $\ww = (w_1, \ldots, w_n)$ as
the {\em weight
vector\/} of the chain.
For extensive treatment of MTF, see \cite{fill-mtf},
which contains
pointers to the sizable literature on the subject.  Hendricks
\cite{hendricks} was the first
to show that the stationary distribution of the MTF Markov chain is
given by~(\ref{target}).

Fill \cite{fill-mtf-ps} used MTF as a case study to compare CFTP and
FMMR.  We omit many
details but for completeness describe the set-up briefly.  Partially
order the symmetric
group ${\SSS}_n$ by declaring $z \leq z'$ if $z'$ can be
obtained from $z$ by a
sequence of adjacent transpositions which switch records out of order
(that
begin in natural
order).  This is the {\em weak Bruhat order\/}.  With
$$ \zh := \mbox{id} = (1,2, \ldots, n), ~~~~ \oh := \mbox{rev}
= (n, n-1, \ldots, 1),$$
we have $\zh \leq z \leq \oh$ for all $z \in {\SSS}_n$.  (For the
definition of the
{\em Bruhat order\/}, used later, delete the word ``adjacent.'')  The
MTF chain
possesses the following monotone transition rule with respect to the
weak Bruhat order.
Let~$U$ be a random variable satisfying
$P(U=i) = w_i$ for $1 \leq i \leq n$.  Define
$$\phi(z,i) = \mbox{move}_i(z) \mbox{ for } z \in {\SSS}_n \mbox{
and } 1 \leq i
\leq n, $$
where $\mbox{move}_i(z)$ is defined to be the permutation resulting from
the list~$z$ by
requesting record~$i$ and applying the MTF rule.  It is easily checked
(see Lemma~2.2
in~\cite{fill-mtf-ps}) that this gives a monotone transition rule for~$M$.

MTF, of course, is not a reversible Markov chain; however,
it is relatively straightforward to generate transitions from the
time-reversed chain.
We refer the reader to~\cite{fill-mtf-ps} for further details
on implementing MTF both using CFTP and using FMMR.

Our first result (Theorem~\ref{revisit})
exhibits explicitly the dependence of~$\Tf$ on the initial
state~$\Zf$.  In what follows, given $z \in \SSS_n$,
let $y_r := w_{z_r}$ for $1 \leq r \leq n$.  In this notation,
(\ref{target}) can be written in the form
$$
\pi(z) = \prod_{r = 1}^n \frac{y_r}{1 - y^+_{r - 1}},
$$
where we have also introduced the notation
$$
y^+_r := \sum_{j = 1}^r y_j,\ \ 0 \leq r \leq n,
$$
for any vector $(y_1, \ldots, y_n)$.

\begin{theorem}\label{revisit}

{\rm (a)}~The conditional distribution of~$\Lc(\Tf)$ given the
initial state~$\Zf$ satisfies
$$
\Lc(\Tf\,|\,\Zf = z) = \Lc(T_{z}),
$$
where the law of~$T_{z}$ is the convolution of
{\rm Geometric}$(1 - y^+_r)$ distributions,
$0 \leq r \leq n - 2$.  We write
$$
T_{z} \sim \oplus_{r = 0}^{n - 2} \mbox{\rm Geom$(1 - y^+_r)$}.
$$

{\rm (b)}~The random variables~$T_{z}$ decrease stochastically in the
Bruhat order for~$z$.

{\rm (c)}~The distribution $\Lc(T_{z})$ is stochastically minimized
(respectively, maximized)
by choosing $z = \rev$ (resp.,\ $z = \id$).  In that case we find
\begin{equation}
\label{Trevid}
T_{\rev} \sim \oplus_{r = 2}^{n}     \mbox{\rm Geom$(w^+_r)$\ \ [resp.,\ }
T_{\id}  \sim \oplus_{r = 0}^{n - 2} \mbox{\rm Geom$(1 - w^+_r)$].}
\end{equation}
\end{theorem}

\begin{proof}
Part~(a) is a consequence of~(\ref{run2}) in our Theorem~\ref{runtime}
and Lemma~3.7 in~\cite{fill-mtf-ps}; indeed, that lemma states that
$\Lc(\Tc\,|\,\Zc = z) =
\oplus_{r=0}^{n-2} \mbox{Geom}(1 - y_r^+)$.  For the {\em weak\/}
Bruhat order, Lemma~3.9
in~\cite{fill-mtf-ps} gives part~(b); but one need only compare $\Lc(T_{z})$ and
$\Lc(T_{z'})$ when~$z$ and~$z'$ differ by any transposition to see that part~(b)
holds for the Bruhat order.  Part~(c) is an immediate consequence of
part~(b).
\end{proof}

\begin{remark}
Theorem~\ref{revisit} for the special case of MTF belies the general
Remark~\ref{nogeneral}.
Indeed, for {\em every\/} initial state for FMMR, we know exactly the
distribution of~$\Tf$.
The {\em worst\/} starting state for coalescence is the identity
permutation, and the result
for $\Lc(T_{\id})$ in
Theorem~\ref{revisit}(c) recaptures Theorem~4.2 in~\cite{fill-mtf-ps}.
The comparison of FMMR and CFTP in~\cite{fill-mtf-ps} was based on
starting FMMR in this worst
state.  In the next section we will discuss how much speedup can be
achieved by instead starting
in the {\em best\/} state, the permutation~$\rev$.
\end{remark}

\sect{Speedup results for MTF}\label{speedup}

\subsection{General weight vectors} \label{general}

From now on, we abbreviate~$T_{\rev}$ of~(\ref{Trevid}) as~$T$.
We first consider how~$\Lc(T)$ varies with the weight vector~$\ww$.
For an understanding of the terminology used in Theorem~\ref{global}
and a thorough treatment
of majorization, see~\cite{MandO}.

\begin{theorem}\label{global}
The distribution $\oplus_{r = 2}^{n} \mbox{\rm Geom$(w^+_r)$}$ of~$T$ is
stochastically strictly Schur-concave in the weight vector~$\ww$.
In particular, the distribution is stochastically maximized
(respectively, minimized),
over all vectors~$\ww$ with $w_1 \geq w_2 \geq \cdots \geq w_n \geq 0$, at
the uniform distribution $\ww = (1/n, \ldots, 1/n)$ (resp.,\ at any
distribution~$\ww$ with $w_1 + w_2 = 1$).
\end{theorem}

\begin{proof}
The result follows simply from the fact that the Geometric$(p)$
distribution is stochastically strictly decreasing in~$p$.
\end{proof}

\begin{remark}
\label{extremes}

(a)~The possibility $w_1 + w_2 = 1$ is ruled out for MTF by our
assumption
that all weights are positive.  Nevertheless it is a limiting case.
In this limiting case, $T = n - 1$ with probability~$1$.
At the other extreme of uniform weights, asymptotics for~$\Lc(T)$ are
well known (since this is a
slight modification of the standard coupon collector's problem).
(The distribution of~$\Lc(T)$ for uniform weights is treated in detail
in Theorem~4.3(a) and Section~4.2 of~\cite{fill-mtf}.
Roughly put, the distribution of~$T$ is concentrated tightly about
$n \ln n$.)  Thus,
for {\em any\/} sequence $\ww^{(1)} = (w_{11}), \ww^{(2)} =
(w_{21}, w_{22}, \ldots), \ldots$ of
weight vectors, writing $T \equiv T_n$ for the~$T$ corresponding to
weight vector~$\ww^{(n)}$
we have 
$$
P(T \geq n - 1) = 0 ~~~~ \mbox{ and } ~~~~
\liminf_{c \to - \infty}\,\liminf_{n \to \infty} P(T \leq n \ln n + c  n) = 1.
$$
So the distribution of~$T$ is always tightly sandwiched between~$n$
and about $n \ln n$,
in sharp contrast (cf.\ Table~1 of~\cite{fill-mtf-ps} or Table~1 below) to the
distribution of~$T_{\id}$ or of~$\Tc$.

(b)~According to Remark~2.6 and the sentence following~(3.2)
in~\cite{fill-mtf-ps}, $\Lc(\Tc)$
is strictly Schur-{\em convex\/} in~$\ww$.  In particular, the
{\em best-case\/} $\Lc(\Tc)$,
corresponding to uniform weights, equals the {\em worst-case\/}
$\Lc(T)$, also corresponding to
uniform weights.

\end{remark}

\subsection{Specific examples of weight vectors} \label{specific}

We
now measure quantitatively, for certain standard examples of weight
vectors, the speedup gained
for FMMR by using the best choice of initial permutation, $\rev$.
Given a triangular array of weights $\ww^{(n)} = (w_{ni}, i = 1,
\ldots, n)$, $n \geq 1$,
we say that $k_n$ steps are necessary and sufficient for convergence
of~$\Lc(T)$ to mean that 
$$
\frac{T_n}{k_n} \to 1 ~~ \mbox{ in probability.}
$$
Here~$T_n$ denotes~$T_{\rev}$ for the weight vector $\ww^{(n)}$;
when there is no danger of confusion, we will sometimes drop
the subscript~$n$.

For some examples of $\ww^{(n)}$ we can obtain results of sharper form
than provided by ``$k_n$ steps are necessary and sufficient''.
However, for simplicity and for uniformity of presentation, we stick
to the above definition.

Let $H_n^{(\alpha)} := \sum_{i = 1}^n i^{-\alpha}$ for $\alpha > 0$,
and let $\zeta(\alpha) := \sum_{i=1}^{\infty} i^{-\alpha}$, $\alpha > 1$,
denote the Riemann zeta function.
We consider the following choices of
weights, where (now suppressing dependence on~$n$ in our notation)
each weight vector~$\ww$ is listed up to a constant of
proportionality.
The numbers of steps necessary and sufficient for convergence
of~$\Lc(T)$
for these examples are stated in Theorem~\ref{specificthm} and
collected in Table~1.  The second and third columns of Table~1 are taken  
from \cite{fill-mtf-ps}.
[In these columns, the meaning of ``$c k_n$ steps are necessary and sufficient''
is that, for some~$h$ and~$H$,
$$
h(c) \leq \liminf_{n \to \infty} P(T_n \leq \lf c k_n \rf) \leq \limsup_{n \to \infty} P(T_n \leq
\lf c k_n \rf) \leq H(c),
$$
where $0 < h(c) \leq H(c) < 1$ for all $c \in (0, \infty)$, $h(c) \to  0$ as $c \to 0$, and $H(c)
\to 1$ as $c \to \infty$.]
The fourth column
in Table~1 is the content of our next theorem.

\vspace{.1in}
\begin{center}
\begin{tabular}{||l|l||} \hline
{\rm Weights} & {\rm $w_i \propto$} \\ \hline
{\rm Uniform} & $1$                 \\
{\rm Zipf's law}
               & $i^{-1}$            \\
{\rm Generalized Zipf's law (GZL)}
               & {\rm $i^{-\alpha}$,\ \ $\alpha > 0$ fixed}
                                     \\
{\rm Power law}
               & {\rm  $(n-i+1)^s$,\ \ $s > 0$ fixed}
                                     \\
{\rm Geometric}
               & {\rm $\theta^i$,\ \ $0 < \theta < 1$ fixed}
                                     \\
\hline
\end{tabular}
\end{center}

\vspace{.1in}
\begin{center}
\small Table 1.\ Rates of convergence for~$\Lc(\Tc)$ and~$\Lc(\Tf)$. \\
\ \\
\normalsize
\begin{tabular}{||l|l|l||l|} \hline
{\rm Weights} & {\rm  $\Lc(\Tc)$}
                                   & {\rm $\Lc(\Tf) (worst)$} & {\rm $\Lc(\Tf) (best)$}          
\\ \hline {\rm Uniform} & $n \ln n $   & $n \ln n $ & $n \ln n$                \\
{\rm Zipf's law}
               & $n (\ln n)^2$
                                   & $n (\ln n)^2$  & $n$     \\
{\rm GZL}
               &                   &                       &         \\
~  $0 < \alpha < 1$
               & $\frac{n}{1 - \alpha} \ln n $
                                   & $\frac{n}{1 - \alpha} \ln n $ & $\frac{n}{\alpha}$
                                                                    \\
~  $\alpha > 1$
               & $\zeta(\alpha) n^{\alpha} \ln n$
                                   & $\zeta(\alpha) n^{\alpha} \ln n$ & $n$
                                                                    \\
{\rm Power law}
               & $c n^{s+1}$       & $c n^{s+1}$      & $\frac{n \ln n}{s+1}$              \\
{\rm Geometric}
               & $c \theta^{-n}$   & $c \theta^{-n}$    & $n$       \\
\hline
\end{tabular}
\end{center}
\vspace{.1in}

\begin{theorem}
\label{specificthm}
\
{\rm (a) (Uniform weights.)} If $w_i = 1/n$ for all $i$, then
$n \ln n$ steps are necessary and sufficient.

{\rm (b) (Zipf's law.)} If $w_i = (H_n i)^{-1}$, with 
$H_n := H_n^{(1)} = \sum_{k=1}^n k^{-1}$,
then $n$ steps are necessary and sufficient.

{\rm (c) (Generalized Zipf's law.)} When $w_i = \left( i^{\alpha} H_n^{(\alpha)}
\right)^{-1}$, {\rm (i)} if $0 < \alpha < 1$, then $n/\alpha$ steps are necessary and
sufficient, and {\rm (ii)} if $\alpha > 1$, then $n$ steps are necessary and sufficient.

{\rm (d) (Power law.)} Fix $s > 0$.  If $w_i = (n-i+1)^s/f(n,s),$
with $f(n,s) := \sum_{j=1}^n j^s$, then $ \frac{n \ln n}{s + 1}$ steps are necessary and
sufficient.

{\rm (e) (Geometric  weights.)} Fix $0 < \theta < 1$. If
$w_i = (1 - \theta) \theta^{i - 1}$ for $i = 1, \ldots, n-1$ and
$w_n = \theta^{n - 1}$,
then $n$
steps are necessary and sufficient.
\end{theorem}

\begin{proof}

We shall ignore the trivialities induced by the need to consider integer
parts in various arguments,
leaving these to the meticulous reader.
\smallskip

(a)~(Uniform weights.)
The asymptotics here are well known, as this is
essentially the standard coupon collector's problem.
A very sharp asymptotic result is that
$$
P(T >\lfloor n \ln n + c n \rfloor) \to 1 - (1+e^{-c}) e^{-e^{-c}}, \qquad c \in \RR.
$$
A thorough treatment of the uniform-weights case is provided by Diaconis et al.~\cite{DFP}.

We establish the remaining results
[as we could also have established~(a)]
by showing, in each case,
(i) that
the number of steps $k_n$
claimed to be necessary and sufficient is the
lead order term
of $\EE[T]$, that is,
that $\EE[T_n] \sim k_n$
as $n \to \infty$, 
and~(ii) that the standard deviation
of~$T_n$ is $o(\EE[T_n])$.
The result then follows
by application of
Chebyshev's inequality.
Showing~(ii) for each
of the weight examples is easy since
\begin{eqnarray*}
\Var[T] & = & \sum_{r=2}^n \frac{1 - w_r^+}{(w_r^+)^2}
  = \sum_{r=2}^n \frac{1}{(w_r^+)^2} - \EE[T] \\
& \leq & \frac{1}{w_2^+} \sum_{r=2}^n \frac{1}{w_r^+} - \EE[T]
  = \EE[T]\left( \frac{1}{w_2^+} - 1 \right),
\end{eqnarray*}
and it is easy to check in each case that $1 / w_2^+ = o(\EE[T])$.
The remainder of the proof
thus consists of showing~(i). In each
case we give explicit upper and 
lower bounds for $\EE[T]$.

(b)~(Zipf's law weights.) Here
\begin{eqnarray}
n - 1 \leq \EE[T]
&   =  & H_n \sum_{r=2}^n (H_r)^{-1} \leq (\ln n + 1)
            \sum_{r=2}^n \frac{1}{\ln(r+1)} \nonumber \\
\label{Zipf1}
& \leq & (\ln n + 1) \int_{2}^{n + 1}\!\frac{dx}{\ln x}.
\end{eqnarray} 
Observe that
$$
\int_{2}^{n + 1}\!\frac{dx}{\ln x}
  = \frac{n + 1}{\ln(n + 1)} - \frac{2}{\ln  2}
      + \int_{2}^{n + 1}\!\frac{dx}{(\ln x)^2},
$$
and that
$$
\int_2^{n / (\ln n)^2}\!\frac{dx}{(\ln x)^2} \leq \frac{n / (\ln n)^2}{(\ln 2)^2}
$$
and
$$
\int_{n / (\ln n)^2}^{n + 1}\!\frac{dx}{(\ln x)^2}
\leq \frac{1}{\ln [n / (\ln n)^2]} \int_{n / (\ln n)^2}^{n + 1}\!\frac{dx}{\ln x}
\leq \frac{1}{\ln[n / (\ln n)^2]} \int_2^{n + 1}\!\frac{dx}{\ln x}.
$$
Thus, 
$$
\int_2^{n + 1}\!\frac{dx}{\ln x}
\leq \frac{n + 1}{\ln(n + 1)} 
  - \frac{2}{\ln 2} + \frac{n}{(\ln n)^2 (\ln 2)^2}
+ \frac{1}{\ln[n / (\ln n)^2]} \int_2^{n + 1}\!\frac{dx}{\ln x},
$$
i.e.,
\begin{eqnarray*}
\int_2^{n + 1}\!\frac{dx}{\ln x}
  &\leq& \left[ 1 - \frac{1}{\ln[n / (\ln n)^2]} \right]^{-1}
           \left[ \frac{n + 1}{\ln(n + 1)}
           + \frac{n}{(\ln n)^2 (\ln 2)^2}  - \frac{2}{\ln 2} \right]\\
  &  = & \frac{n + 1}{\ln(n + 1)} + O \left( \frac{n}{(\ln n)^2} \right).  
\end{eqnarray*}
Continuing now from~(\ref{Zipf1}) we have
$$
\EE[T] 
\leq (\ln n + 1) \left( \frac{n + 1}{\ln(n + 1)} \right)
+ O \left( \frac{n}{\ln n} \right)
= n + O \left( \frac{n}{\ln n} \right).
$$

(c)~(Generalized Zipf's law.)
For $0 < \alpha < 1$,
\begin{eqnarray*}
\EE[T] & = & H_n^{(\alpha)} \sum_{r=2}^n \frac{1}{H_r^{(\alpha)}}
\leq (n - 1)^{1-\alpha} \sum_{r=2}^n
\frac{1}{(r+1)^{1-\alpha} -1} \\
& \leq & n^{1-\alpha} \sum_{r=3}^{n + 1} \frac{1}{r^{1-\alpha} -1}
= n^{1-\alpha} \sum_{r=3}^{n + 1} \frac{1}{r^{1-\alpha}} 
\left( 1 + O\left( r^{\alpha-1} \right) \right) \\
& = & n^{1-\alpha} \sum_{r=3}^{n+1} \left( r^{\alpha -1} + O( r^{2\alpha-2}) \right)
= n^{1-\alpha}\sum_{r=3}^{n + 1} r^{\alpha-1} + O(c_n) \\
& = & \frac{n}{\alpha} + o(n),
\end{eqnarray*}
where $c_n$ is defined as $n^{1 - \alpha}$ if $0 < \alpha < 1 / 2$, as $n^{1 / 2} \ln n$ if $\alpha = 1 /
2$, and as $n^{- \alpha}$ if $1 / 2 < \alpha < 1$.
Also,
\begin{eqnarray*}
\EE[T] & \geq & \left( (n+1)^{1 - \alpha} -1 \right) \sum_{r=2}^n \frac{1}{(r - 1)^{1-\alpha}} \\
& = & (1 + o(1)) n^{1 - \alpha} \frac{n^{\alpha}}{\alpha} \\ 
& = & \frac{n}{\alpha} + o(n).
\end{eqnarray*}

For $\alpha > 1$, we use the fact that
\begin{equation}
\label{zeta}
  H_n^{(\alpha)} = \zeta(\alpha) - \frac{n^{-(\alpha -1)}}{\alpha -1} + O(n^{-\alpha}).
\end{equation}
Now
\begin{eqnarray*}
\sum_{r=2}^n \frac{1}{H_r^{(\alpha)}} & = & \sum_{r=2}^n \left[
\zeta(\alpha) - \frac{r^{-(\alpha -1)}}{\alpha -1} + O(r^{-\alpha}) \right]^{-1} \\
& = & \frac{1}{\zeta(\alpha)} \sum_{r=2}^n \left[
1 - \frac{r^{-(\alpha -1)}}{(\alpha -1)\zeta(\alpha)} + O(r^{-\alpha}) \right]^{-1} \\
& = & \frac{1}{\zeta(\alpha)} \sum_{r=2}^n \left[
1 + \frac{r^{-(\alpha -1)}}{(\alpha-1)\zeta(\alpha)} + O(r^{-\alpha})
+ O(r^{-(2\alpha - 2)}) \right] \\
& = & \frac{n}{\zeta(\alpha)} + o(n).
\end{eqnarray*}
Together with (\ref{zeta}) this gives
\begin{eqnarray*}
\EE[T] & = & H_n^{(\alpha)} \sum_{r=2}^n \frac{1}{H_r^{(\alpha)}}  \\
& = & \left[ \zeta(\alpha) + O(n^{-(\alpha-1) } ) \right] \left[\frac{n}{\zeta(\alpha)}
+ o(n) \right] \\
& = & n + o(n).
\end{eqnarray*}

\vspace{.1in}

(d)~(Power law.)  Here 
\begin{eqnarray*}
\EE[T]
&  = & \sum_{r=2}^n \frac{f(n,s)}{(n-r+1)^s + \cdots + n^s}
  \leq  f(n,s) \sum_{r=2}^n \frac{s+1}{n^{s+1} - (n-r)^{s+1}} \\
&  = & \frac{f(n,s) (s+1)}{n^{s+1}} \sum_{r=2}^n \frac{1}{1 - (1 -
          \frac{r}{n})^{s+1}}.
\end{eqnarray*}
The 
inequality follows from an integral comparison.  Another
integral comparison shows that the last 
sum above is bounded above by 
$$
\int_1^n\!\frac{dx}{1-(1-\frac{x}{n})^{s+1}}
= n \int_0^{1 - \frac{1}{n}}\!\frac{dy}{1- y^{s+1}} =: n \times I.
$$
Now 
\begin{eqnarray*}
I
&  = & \int_0^{1 - \frac{1}{n}} \sum_{k=0}^{\infty} (y^{s+1})^k\,dy \\
&  = & \sum_{k=0}^{\infty} \frac{1}{(s+1)k + 1} \left(1- \frac{1}{n} \right)^{(s+1)k + 1} \\
&\leq& \left( 1- \frac{1}{n} \right) + \frac{1-\frac{1}{n}}{s+1}
          \sum_{k=1}^{\infty}
          \frac{ \left[  (1- \frac{1}{n})^{s+1}  \right]^k }{k} \\
&  = &\left( 1- \frac{1}{n} \right) + \frac{1-\frac{1}{n}}{s+1}
          \left| \ln \left( 1- \left( 1- \frac{1}{n} \right)^{s+1} \right) \right| \\
&\leq& \left( 1- \frac{1}{n} \right)
          + \frac{1-\frac{1}{n}}{s+1}
          \left| \ln \left( \frac{s+1}{n} - \frac{(s+1)s}{2n^2} \right) \right|\\
&\leq& \left( 1- \frac{1}{n} \right)
          + \frac{1-\frac{1}{n}}{s+1} \ln n \\
& \leq & \frac{\ln n}{s+1} + 1,
\end{eqnarray*}
where the penultimate inequality holds for all sufficiently large~$n$  [in particular,
for $n \geq (s + 1) / 2$].
We then have that 
\begin{eqnarray*}
\EE[T] & \leq & \frac{f(n,s) (s+1)}{n^s}
\left( \frac{\ln n}{s+1} + 1 \right) \\
& \leq & \frac{n \ln n}{s+1} + n + \ln n + s+1 = \frac{n \ln n}{s+1} + O(n) = (1 + o(1)) \frac{n
\ln n}{s + 1},
\end{eqnarray*}
using 
$$
f(n,s) \leq \frac{n^{s+1}}{s+1} + n^s.
$$

For the lower bound,
\begin{eqnarray*}
\EE[T] &\geq& f(n, s) \sum_{r = 2}^n
                 \frac{s + 1}{(n + 1)^{s + 1} - (n + 1 - r)^{s + 1}} \\
        &  = & \frac{f(n, s) (s + 1)}{(n + 1)^{s + 1}} \sum_{r = 2}^n
                 \left[ 1 - \left( 1 - \frac{r}{n + 1} \right)^{s + 1} \right]^{-1}.
\end{eqnarray*}
But
\begin{eqnarray*}
\sum_{r = 2}^n \left[ 1 - \left( 1 - \frac{r}{n + 1} \right)^{s +  1} \right]^{-1}
&\geq& \int_2^{n + 1} \left[ 1 - \left( 1 - \frac{x}{n + 1} \right)^{s + 1} \right]^{-1}\,dx \\
&  = & (n + 1) \int_0^{1 - \frac{2}{n + 1}} \frac{dy}{1 - y^{s + 1}}.
\end{eqnarray*}
We can show, but omit the details, that
$$ \int_0^{1 - \frac{2}{n + 1}} \frac{dy}{1 - y^{s + 1}} = \frac{\ln n}{s + 1} - O(1),$$
%
%
so
$$
\EE[T] \geq \left( \frac{n}{n + 1} \right)^{s + 1} \frac{(n + 1) \ln n}{s + 1} - O(n) \geq
\frac{n \ln n}{s + 1} - O(n),
$$
where we have used
$$
f(n, s) \geq \frac{n^{s+1}}{s+1}.
$$

(e)~(Geometric weights.) We have
\begin{eqnarray*}
n - 1 \leq \EE[T]
&    = & n - 1 + \sum_{r = 2}^{n - 1} \frac{\theta^r}{1 - \theta^r} \\
& \leq & n - 1 + \sum_{r = 2}^{\infty} \frac{\theta^r}{1 - \theta^r} \\
&    = & n + O(1) = (1 + o(1)) n.
\end{eqnarray*}
\end{proof}

\section{Coalescence into a set}\label{last}
Here we present a different approach to speeding up the FMMR
algorithm,
which has the same spirit as the other results in this paper.  In brief,
recall that FMMR  starts in a user-chosen state and then, in the
second phase of the algorithm,  checks whether there is coalescence
to that state.  
In the generalization we consider here, one starts
the algorithm in some subset of the state space (not necessarily a
singleton) and then checks if there is coalescence back to that set
(but not necessarily to the state in which the algorithm began).

\begin{theorem}
\label{coalS}
In the general setting for FMMR described in Section~\ref{subfmmr},
let~$S$ be a subset of the state space
and define $\pi_0(\cdot) := \pi(\cdot | S)$.  Consider the modified algorithm which
starts in a state~$\XX_0$ distributed according to~$\pi_0$
and outputs $W := \XX_{- T}$,
where~$T$ is defined to be the smallest~$t$ such that all the forward
trajectories from time $- t$ coalesce into~$S$, \ie, such that
$\YY^{(- t)}(x) \in S$ for every state~$x$.  Then~$W$ has the stationary
distribution~$\pi$.  Further, the algorithm is interruptible
(\ie, $T$ and~$W$ are independent random variables).
\end{theorem}

\begin{proof}
For simplicity we consider only the discrete case.
It suffices to show that
\begin{equation}
\label{suffices}
P(T \leq t, \XX_{- t} = x) = P(T \leq t) \pi(x)
\end{equation}
for every~$t$ and~$x$, for then
\begin{eqnarray*}
P(T = t,\ W = x)
 &=& P(T = t,\ \XX_{- T} = x) = P(T = t,\ \XX_{- t} = x) \\
 &=& P(T \leq t,\ \XX_{- t} = x) - P(T \leq t - 1,\ \XX_{- t} = x) \\
 &=& P(T \leq t) \pi(x) - P(T \leq t - 1) \pi(x) \\
 &=& P(T = t) \pi(x),
\end{eqnarray*}
as desired.  Here we have used the fact that~$\pi$ is stationary for
the time-reversed kernel~$\Kt$, so that
\begin{eqnarray*}
\lefteqn{P(T \leq t - 1,\ \XX_{- t} = x)} \\
 &=& \sum_y P(T \leq t - 1, \XX_{- (t - 1)} = y)\,P(\XX_{- t} = x\,|\,T
       \leq t - 1,\ \XX_{- (t - 1)} = y) \\
 &=& \sum_y P(T \leq t - 1) \pi(y) \Kt(y, x)\mbox{\quad by~\eqref{suffices}
       and the Markov property for~$\Kt$} \\
 &=& P(T \leq t - 1) \pi(x).
\end{eqnarray*}

To establish~\eqref{suffices}, we first observe that
\begin{eqnarray}
P(\XX_{- t} = x) 
 &=& \sum_z \pi_0(z) \Kt^t(z, x)
  =  \pi(x) \sum_z \frac{\pi_0(z)}{\pi(z)} K^t(x, z) \nonumber \\
\label{uncond}
 &=& \frac{\pi(x)}{\pi(S)} \sum_{z \in S} K^t (x, z)
  =  \frac{\pi(x)}{\pi(S)} K^t (x, S).
\end{eqnarray}
One can check that, conditionally given $\XX_{- t} = x$,
the forward trajectory $(\XX_{- t}, \ldots,
\XX_0)$ has the same distribution as a $K$-trajectory
conditioned to start at~$x$ and end in~$S$. 
Therefore, by the algorithm's design,
\begin{eqnarray*}
\lefteqn{P(T \leq t\,|\,\XX_{-t} = x)} \\
 &=& P(\mbox{forward coalescence into~$S$ over a time-interval
       of length~$t$}) / K^t(x, S).
\end{eqnarray*}
Combining this with~\eqref{uncond} we
conclude~\eqref{suffices}, and the additional result
$$
P(T \leq t) = P(\mbox{forward coalescence into~$S$ over
                a time-interval of length~$t$}) / \pi(S).
$$
\end{proof}

\begin{remark}
In the monotone case, if $S$ is a down-set
(meaning:\ $z \in S$ and $y \leq z$ implies $y \in S$),
then the computational problem of determining whether or not
there is coalescence into~$S$ is eased considerably:\ we need
only determine whether the terminal state (call it~$y$) of
the forward trajectory started in~$\oh$ belongs to~$S$.
And so if~$S$ is a \emph{principal} down-set, that is, if 
$S = \{z: z \leq z_0\}$ for some~$z_0$, the problem is even easier:\ 
we need only check whether $y \leq z_0$.
\end{remark}

We will now give a ``toy'' application of these ideas to MTF
by describing an algorithm to build up a stationary observation in
just $n - 1$ steps, regardless of the weights $w_1, \ldots, w_n$.
Let $\pi_k$ denote the MTF stationary distribution
on ${\cal S}_k$ 
restricted to the (normalized) weights $w_1, \ldots, w_k$; that is,
to the weights
$w_1/w_k^+, \ldots, w_k/w_k^+$. Let $\mbox{MTF}_k$ denote the MTF
process on ${\cal S}_k$,
and let $S_k$
denote the set of permutations of $\{1, \ldots, k\}$ that begin with~$k$.
Observe that $S_k$ is the
principal order ideal $\{z: z \leq z_k\}$
in the dual (\ie, ``upside-down'') Bruhat order
of ${\cal S}_k$,
where $z_k$ is the permutation $(k\ 1\ \cdots\ k - 1)$.
(We will not refer to the symmetric group ${\cal S}_k$ any further; thus there will be
no notational confusion with its special subset $S_k$.)
Inductively, after $k$ steps of our algorithm we will have a permutation
distributed according to $\pi_{k + 1}$; thus, after $n - 1$ steps we will have
an observation from~$\pi$.
 
Initialize the algorithm (step 0) with the permutation $(1)$ on $\{1\}$. 
Suppose that after $k-1$ steps we have the permutation $x = (x_1, \ldots,  x_k)$
distributed according to $\pi_k$.  For the next ($k$th) step,
we first get immediately an observation from
$\pi_{k+1}(\cdot\,|\,S_{k+1})$, namely, $(k+1, x_1, \ldots, x_k)$.  Then we apply
the ``coalesce into~$S$'' routine of \refT{coalS}, taking $S = S_{k + 1}$.
We claim that that routine will terminate in a single step!
Indeed, in one time-reversed $\mbox{MTF}_{k+1}$ transition we obtain
the  permutation
$$
x' = (x_1, \ldots, x_{j-1}, k+1, x_j, \ldots, x_k)
$$
for some $j$.  In the forward phase of the routine, record $k+1$ is brought to the front
of
every trajectory,
giving coalescence into the set $S_{k+1}$.  We thus conclude from \refT{coalS}
that $x' \sim \pi_{k+1}$, completing the induction. 

\sect{Acknowledgements}

This research was carried out in part while the first-listed author
was a member of the Department of Mathematics and Computer Science
at Clarkson University, and while the second-listed author was Visiting Researcher,
Theory Group, Microsoft Research.

\end{document}